\documentclass{article}[12pt]

\usepackage{latexsym}
\usepackage{amsmath}
\usepackage{amssymb}

\begin{document}

\title{EXPLICITLY SOLVABLE SYSTEMS OF FIRST-ORDER ORDINARY DIFFERENTIAL
EQUATIONS WITH POLYNOMIAL RIGHT-HAND SIDES, AND THEIR PERIODIC VARIANTS}

\author{Francesco Calogero$^{a,b}$\thanks{e-mail: francesco.calogero@roma1.infn.it}
\thanks{e-mail: francesco.calogero@uniroma1.it}
 , Farrin Payandeh$^c$\thanks{e-mail: farrinpayandeh@yahoo.com}
 \thanks{e-mail: f$\_$payandeh@pnu.ac.ir}}

\maketitle   \centerline{\it $^{a}$Physics Department, University of
Rome "La Sapienza", Rome, Italy}

\maketitle   \centerline{\it $^{b}$INFN, Sezione di Roma 1}

\maketitle

\maketitle   \centerline{\it $^{c}$Department of Physics, Payame
Noor University, PO BOX 19395-3697 Tehran, Iran}

\maketitle

\begin{abstract}

In this Letter we identify \textit{special} systems of (an \textit{arbitrary}
number) $N$ of first-order Ordinary Differential Equations with \textit{%
homogeneous} polynomials of \textit{arbitrary} degree $M$ on their
right-hand sides, which feature \textit{very simple explicit }solutions; as
well as variants of these systems---with right-hand sides no more
homogeneous---which feature \textit{periodic} solutions. A novelty of these
findings is to consider \textit{special} systems characterized by \textit{%
constraints} involving both their parameters and their initial data.

\end{abstract}

The general system of an \textit{arbitrary} number $N$ of first-order
Ordinary Differential Equations (ODEs) with \textit{homogeneous} \textit{%
polynomials} of \textit{arbitrary} degree $M$ on their right-hand sides
reads as follows:
\begin{subequations}
\label{1}
\begin{eqnarray}
\dot{z}_{n}\left( t\right) &=&\sum_{m_{\ell }}{}^{\left( M\right) }\left\{
c_{nm_{1}m_{2}\cdot \cdot \cdot m_{N}}\left[ z_{1}\left( t\right) \right]
^{m_{1}}\left[ z_{2}\left( t\right) \right] ^{m_{2}}\cdot \cdot \cdot \left[
z_{N}\left( t\right) \right] ^{m_{N}}\right\} ~,  \nonumber \\
n &=&1,2,...,N~,  \label{zndot}
\end{eqnarray}%
where (above and below) the symbol $\sum_{m_{\ell }}{}^{\left( M\right) }$
denotes a sum running over \textit{all nonnegative} values of the $N$
indices $m_{\ell }$ subject to the restriction%
\begin{equation}
\sum_{\ell =1}^{N}\left( m_{\ell }\right) =M~,  \label{M}
\end{equation}%
implying that the polynomials in $N$ variables $z_{n}\left( t\right) $ in
the right-hand sides of the $N$ ODEs (\ref{zndot}) are \textit{all
homogeneous} of degree $M$.

\textbf{Notation}. Throughout this paper $M$ and $N$ are \textit{positive}
integers larger than \textit{unity}; the index $n$ takes \textit{positive}
\textit{integer} values; indices and exponents such as $m_{1}$, $m_{2}$,...
take \textit{all} the \textit{nonnegative integer} values consistent with
the restriction (\ref{M}); the independent variable $t$ can be considered as
playing the role of "time", taking \textit{all nonnegative real} values (but
it shall also be eventually convenient to replace it formally with the
\textit{complex} variable $\tau $, see below); a superimposed dot indicates
a $t$-differentiation; the coefficients $c_{nm_{1}m_{2}\cdot \cdot \cdot
m_{M}}$ are ($t$-independent) parameters; while of course the dependent
variables $z_{n}\equiv z_{n}\left( t\right) $ are functions of the
independent variable $t$ and ascertaining their $t$-evolution from the set
of $N$ \textit{initial data} $z_{n}\left( 0\right) $ is our main task. The
coefficients $c_{nm_{1}m_{2}\cdot \cdot \cdot m_{M}}$ and the dependent
variables $z_{n}\left( t\right) $ might be restricted to be \textit{real};
but in the last part of this paper we shall assume that they are instead
\textit{complex}, setting
\end{subequations}
\begin{equation}
c_{nm_{1}m_{2}\cdot \cdot \cdot m_{M}}=a_{nm_{1}m_{2}\cdot \cdot \cdot
m_{M}}+\mathbf{i}b_{nm_{1}m_{2}\cdot \cdot \cdot m_{M}}~;  \label{zclm}
\end{equation}%
and we shall as well replace the independent variable $t$ with a \textit{%
complex} variable $\tau ,$ see below eq. (\ref{xyzn}); here and below of
course $\mathbf{i}$ is the \textit{imaginary unit}, $\mathbf{i}^{2}=-1$.
Finally: below $\omega $ denotes an \textit{arbitrary nonvanishing real}
parameter. $\blacksquare $

The system (\ref{1}) has being investigated over time in an enormous number
of mainly mathematical, or mainly applicative, papers (more than it is
possible to report in an adequate manner: for a seminal paper see, for
instance, \cite{G1960}); although generally for specific, relatively small,
values of $N$ and $M$. The mathematics behind the findings reported in the
present paper is rather elementary; yet these developments may have some
interest---perhaps mainly in applicative contexts---because they are based
on a somewhat \textit{unconventional} approach: to identify \textit{%
explicitly solvable} cases of the system (\ref{1}) by introducing \textit{%
constraints} involving, in addition to the coefficients $c_{nm_{1}m_{2}\cdot
\cdot \cdot m_{M}}$, also the initial data $z_{n}\left( 0\right) $ (which,
in many applicative contexts, may well play the role of \textit{control}
elements, determining the time evolution of the system).

Our main result is the following

\textbf{Proposition}. The system (\ref{1}) features the special solution
\begin{subequations}
\label{Solznt}
\begin{equation}
z_{n}\left( t\right) =z_{n}\left( 0\right) \left( 1+Kt\right) ^{1/\left(
1-M\right) }~,~~~n=1,2,...,N~,  \label{znt}
\end{equation}%
provided there hold the following $N$ explicit \textit{algebraic constraints}
on the \textit{a priori arbitrary }parameter $K$, the coefficients $%
c_{nm_{1}m_{2}\cdot \cdot \cdot m_{M}}$ and the $N$ initial data $%
z_{n}\left( 0\right) $:
\begin{eqnarray}
Kz_{n}\left( 0\right) &=&\left( 1-M\right) \left[ \sum_{m_{\ell }}{}^{\left(
M\right) }\left\{ c_{nm_{1}m_{2}\cdot \cdot \cdot m_{M}}\left[ z_{1}\left(
0\right) \right] ^{m_{1}}\left[ z_{2}\left( 0\right) \right] ^{m_{2}}\cdot
\cdot \cdot \left[ z_{N}\left( 0\right) \right] ^{m_{N}}\right\} \right] ~,
\nonumber \\
n &=&1,2,...,N~.~\blacksquare  \label{Kzn0}
\end{eqnarray}

\textbf{Remark 1}. The proof that (\ref{Solznt}) satisfies the system of
ODEs (\ref{1}) is elementary: just insert (\ref{znt}) in (\ref{zndot}) and
verify that, thanks to (\ref{M}) and (\ref{Kzn0}), the $N$ ODEs (\ref{zndot}%
) are satisfied. $\blacksquare $

\textbf{Remark 2}. The system of $N$ \textit{algebraic} equations (\ref{Kzn0}%
) generally determines---for any given assignment of the \textit{a priori
arbitrary} coefficients $c_{nm_{1}m_{2}\cdot \cdot \cdot m_{M}}$---$N$ out
of the $N+1$ quantities $K$ and $z_{n}\left( 0\right) $; but it is also
possible to select \textit{ad libitum} $N$ elements out of the \textit{%
complete} set of data $K$, $c_{nm_{1}m_{2}\cdot \cdot \cdot m_{M}}$ and $%
z_{n}\left( 0\right) $, and to then consider these \textit{selected}
elements as those to be determined---by the $N$ conditions (\ref{Kzn0})---in
terms of the remaining \textit{arbitrarily assigned} elements in the \textit{%
complete} set of these data. If one chooses to satisfy these $N$ conditions
by solving the $N$ equations (\ref{Kzn0}) for $N$ of the coefficients $%
c_{nm_{1}m_{2}\cdot \cdot \cdot m_{M}}$---or for the parameter $K$ and $N-1$
of the coefficients $c_{nm_{1}m_{2}\cdot \cdot \cdot m_{M}}$---then this
task can be generally performed \textit{explicitly}, since the relevant
\textit{algebraic} equations to be solved are then \textit{linear} in the
unknown quantities; otherwise these determinations require the solution of
\textit{nonlinear }equations, a task which can be performed \textit{%
explicitly} only rarely in an \textit{algebraic} setting; but which can
generally be performed, with \textit{arbitrary} approximation, in a \textit{%
numerical} context. $\blacksquare $

\textbf{Example 1}. Assume for instance $N=2$ and $M=4$, so that the system (%
\ref{1}) reads as follows (note below the notational simplification):
\end{subequations}
\begin{subequations}
\label{Ex1}
\begin{equation}
\dot{z}_{n}\left( t\right) =\sum_{m=0}^{4}c_{nm}\left[ z_{1}\left( t\right) %
\right] ^{4-m}\left[ z_{2}\left( t\right) \right] ^{m}~,~~~n=1,2~,
\end{equation}%
featuring $2$ dependent variables $z_{n}\left( t\right) $ and $10$ \textit{a
priori arbitrary} coefficients $c_{nm}$ ($n=1,2$; $m=0,1,2,3,4$). Then the
solution (\ref{znt}) reads as follows:
\begin{equation}
z_{n}\left( t\right) =z_{n}\left( 0\right) \left[ 1+Kt\right]
^{-1/3}~,~~~n=1,2~,
\end{equation}%
and the $2$ conditions (\ref{Kzn0}) read as follows:
\begin{equation}
Kz_{n}\left( 0\right) =-3\sum_{m=0}^{4}c_{nm}\left[ z_{1}\left( 0\right) %
\right] ^{4-m}\left[ z_{2}\left( 0\right) \right] ^{m}~,~~~n=1,2~.
\end{equation}%
These algebraic constraints can of course be \textit{explicitly} solved for
any $2$ of the $10$ coefficients $c_{nm}$ in terms of the other $8$
coefficients $c_{nm}$ and of the $3$ \textit{arbitrary} data $K$, $%
z_{1}\left( 0\right) $, $z_{2}\left( 0\right) $; or alternatively for $K$
and only $1$ of the $10$ coefficients $c_{nm}$ in terms of the other $9$
coefficients $c_{nm}$ and of the $2$ \textit{arbitrary} initial data $%
z_{1}\left( 0\right) $, $z_{2}\left( 0\right) $; with many other
possibilities left to the imagination of the interested reader. $%
\blacksquare $

The \textit{periodic} variant obtains from the previous results---where we
now assume all quantities to be \textit{complex} and we formally replace the
independent variable $t$ with the \textit{complex} variable $\tau $---via
the following well-known trick (amounting to a simple change of dependent
and independent variables: see, for instance, \cite{C2008}):
\end{subequations}
\begin{equation}
x_{n}\left( t\right) +\mathbf{i}y_{n}\left( t\right) =\left\{ \exp \left[
\mathbf{i}\omega t/\left( M-1\right) \right] \right\} z_{n}\left( \tau
\right) ~,~~~\tau =\left[ \exp \left( \mathbf{i}\omega t\right) -1\right]
/\left( \mathbf{i}\omega \right) ~,  \label{xyzn}
\end{equation}%
implying $\dot{\tau}\left( t\right) =\exp \left( \mathbf{i}\omega t\right) $
and transforming the system (\ref{zndot}) into the following (still \textit{%
autonomous}!) system involving now the $2N$ \textit{real} variables $%
x_{n}\left( t\right) $ and $y_{n}\left( t\right) $ (depending of course on
the \textit{real} independent variable $t$: "time"):
\begin{subequations}
\label{PerSyst}
\begin{eqnarray}
\dot{x}_{n}\left( t\right) &=&-\left[ \omega /\left( M-1\right) \right]
y_{n}\left( t\right) +Re\left[ Z_{n}\left( t\right) \right] ~,  \nonumber
\\
\dot{y}_{n}\left( t\right) &=&\left[ \omega /\left( M-1\right) \right]
x_{n}\left( t\right) +Im\left[ Z_{n}\left( t\right) \right] ~,
\end{eqnarray}%
where (see (\ref{xyzn}) and (\ref{zclm}))
\begin{eqnarray}
Z_{n}\left( t\right) &=&\sum_{m_{\ell }}{}^{\left( M\right) }\left\{ \left(
a_{nm_{1}m_{2}\cdot \cdot \cdot m_{N}}+\mathbf{i}b_{nm_{1}m_{2}\cdot \cdot
\cdot m_{N}}\right) \cdot \right.  \nonumber \\
&&\left. \cdot \left[ x_{1}\left( t\right) +\mathbf{i}y_{1}\left( t\right) %
\right] ^{m_{1}}\cdot \cdot \cdot \left[ x_{N}\left( t\right) +\mathbf{i}%
y_{N}\left( t\right) \right] ^{m_{N}}\right\} ~.
\end{eqnarray}

\textbf{Remark 3}. The fact that \textit{all} solutions $x_{n}\left(
t\right) $, $y_{n}\left( t\right) $ of the system (\ref{PerSyst}) obtained
via the definition (\ref{xyzn}) with $z_{n}\left( \tau \right) $ defined by (%
\ref{znt}) (of course with $t$ replaced there by $\tau $, see (\ref{xyzn}))
are \textit{periodic} with a period $T$ which is an (easily identifiable on
a case-by-case basis) \textit{integer} multiple of the basic period $2\pi
/\left\vert \omega \right\vert $ is rather \textit{obvious}; in case of
\textit{doubt}, see \cite{C2008}. $\blacksquare $

\textbf{Example 2}. As an example of \textit{solvable} system featuring
\textit{periodic} solutions let us display the findings reported in the
\textit{special} case with $N=2$ and $M=4$. Then the system (\ref{PerSyst})
of $4$ ODEs reads as follows:
\end{subequations}
\begin{subequations}
\begin{eqnarray}
\dot{x}_{n}\left( t\right) &=&-\left( \omega /3\right) y_{n}\left( t\right) +%
Re\left[ Z_{n}\left( t\right) \right] ~,~~~n=1,2~,  \nonumber \\
\dot{y}_{n}\left( t\right) &=&\left( \omega /3\right) x_{n}\left( t\right) +%
Im\left[ Z_{n}\left( t\right) \right] ~,~~~n=1,2~,
\end{eqnarray}%
\begin{equation}
Z_{n}\left( t\right) =\sum_{m=0}^{4}\left\{ \left( a_{nm}+\mathbf{i}%
b_{nm}\right) \left[ x_{1}\left( t\right) +\mathbf{i}y_{1}\left( t\right) %
\right] ^{4-m}\left[ x_{1}\left( t\right) +\mathbf{i}y_{1}\left( t\right) %
\right] ^{m}\right\} ~;
\end{equation}%
its \textit{explicit} solutions read as follows:
\end{subequations}
\begin{subequations}
\label{Ex2}
\begin{equation}
x_{n}\left( t\right) =Re\left[ \zeta _{n}\left( t\right) \right]
~,~~~y_{n}\left( t\right) =Im\left[ \zeta _{n}\left( t\right) \right]
~,~~~n=1,2~,
\end{equation}%
\begin{eqnarray}
\zeta _{n}\left( t\right) =\left[ x_{n}\left( 0\right) +\mathbf{i}%
y_{n}\left( 0\right) \right] \exp \left( \mathbf{i}\omega t/3\right) \cdot &&
\nonumber \\
\cdot \left\{ \left[ 1+\left( K_{R}+\mathbf{i}K_{I}\right) \left[ \exp
\left( \mathbf{i}\omega t\right) -1\right] /\left( \mathbf{i}\omega \right) %
\right] ^{-1/3}\right\} ~,~~~n=1,2~, &&
\end{eqnarray}%
provided the $2$ (\textit{a priori arbitrary}) \textit{real} parameters $%
K_{R}$ and $K_{I}$, the $4$ (\textit{a priori arbitrary}) \textit{real}
\textit{initial data} $x_{n}\left( 0\right) $ and $y_{n}\left( 0\right) $
and the $20$ (\textit{a priori arbitrary}) \textit{real} coefficients $%
a_{nm} $ and $b_{nm}$ ($n=1,2;~m=0,1,2,3,4$) are related to each other by
the following $2$ \textit{complex }(i. e., $4$ \textit{real})\textit{\
constraints}:
\begin{eqnarray}
\left( K_{R}+\mathbf{i}K_{I}\right) \left[ \mathbf{x}_{n}\left( 0\right) +%
\mathbf{i}y_{n}\left( 0\right) \right] &&  \nonumber \\
=-3\sum_{m=0}^{4}\left\{ \left( a_{nm}+\mathbf{i}b_{nm}\right) \left[
x_{1}\left( 0\right) +\mathbf{i}y_{1}\left( 0\right) \right] ^{4-m}\left[
x_{1}\left( 0\right) +\mathbf{i}y_{1}\left( 0\right) \right] ^{m}\right\}
~,~n=1,2~.~\blacksquare &&
\end{eqnarray}

\textbf{Final Remark}. As already noted above, the mathematics behind the
results reported above is rather \textit{elementary}. Yet these findings do
not seem to have been advertised so far, while their \textit{applicable}
potential is clearly vast; so---especially among \textit{applied}
mathematicians and \textit{practitioners} of the various scientific
disciplines where systems of ODEs such as those discussed above play a key
role---a wider knowledge of them seems desirable; for instance via their
inclusion in standard compilations of \textit{solvable} ODEs such as \cite%
{PZ}. $\blacksquare $

\textbf{Acknowledgements}. It is a pleasure to thank our colleagues Robert
Conte, Fran\c{c}ois Leyvraz and Andrea Giansanti for very useful
discussions. We also like to acknowledge with thanks $2$ grants,
facilitating our collaboration---mainly developed via e-mail exchanges---by
making it possible for FP to visit twice the Department of Physics of the
University of Rome "La Sapienza": one granted by that University, and one
granted jointly by the Istituto Nazionale di Alta Matematica (INdAM) of that
University and by the International Institute of Theoretical Physics (ICTP)
in Trieste in the framework of the ICTP-INdAM "Research in Pairs" Programme.%
\textbf{\ }Finally, we also\ like to thank Fernanda Lupinacci who, in these
difficult times---with extreme efficiency and kindness---facilitated all the
arrangements necessary for the presence of FP with her family in Rome.

\end{subequations}

\end{document}